\documentclass[11pt]{article}
\usepackage{amsmath,amssymb}

\newtheorem{propo}{{\bf Proposition}}[section]
\newtheorem{coro}[propo]{{\bf Corollary}}
\newtheorem{lemma}[propo]{{\bf Lemma}} \newtheorem{theor}[propo]{{\bf
Theorem}}

\begin{document}

\vspace*{1.0in}

\begin{center} C-SUPPLEMENTED SUBALGEBRAS OF LIE ALGEBRAS 
\end{center}
\bigskip

\begin{center} DAVID A. TOWERS 
\end{center}
\bigskip

\begin{center} Department of Mathematics and Statistics

Lancaster University

Lancaster LA1 4YF

England

d.towers@lancaster.ac.uk 
\end{center}
\bigskip

\begin{abstract}
A subalgebra $B$ of a Lie algebra $L$ is {\em c-supplemented} in $L$ if there is a subalgebra $C$ of $L$ with $L = B + C$ and $B \cap C \leq B_L$, where $B_L$ is the core of $B$ in $L$. This is analogous to the corresponding concept of a  c-supplemented subgroup in a finite group. We say that $L$ is {\em c-supplemented} if every subalgebra of $L$ is c-supplemented in $L$. We give here a complete characterisation of c-supplemented Lie algebras over a general field.
\end{abstract}

\section{Introduction}
\medskip
The concept of a c-supplemented subgroup of a finite group was introduced by Ballester-Bolinches, Wang and Xiuyun in \cite{bwx} and has since been studied by a number of authors. The purpose of this paper is study the corresponding idea for Lie algebras. As we shall see, stronger results can be obtained in this context.
\par
Throughout $L$ will denote a finite-dimensional Lie algebra over a field $F$. 
If $B$ is a subalgebra of $L$ we define $B_L$, the {\em core} (with respect to $L$) of $B$ to be the largest ideal of $L$ contained in $B$. We say that $B$ is {\em core-free} in $L$ if $B_L = 0$. A subalgebra $B$ of $L$ is {\em c-supplemented} in $L$ if there is a subalgebra $C$ of $L$ with $L = B + C$ and $B \cap C \leq B_L$. We say that $L$ is {\em c-supplemented} if every subalgebra of $L$ is c-supplemented in $L$. We shall give a complete characterisation of c-supplemented Lie algebras over a general field. 
\par
Following \cite{gm} we will say that $L$ is {\em completely factorisable} if for every subalgebra $B$ of $L$ there is a subalgebra $C$ such that $L = B + C$ and $B \cap C = 0$. It turns out that c-supplemented Lie algebras are intimately related to the completely factorisable ones, and our results generalise some of those obtained in \cite{gm}. Incidentally, it is claimed in \cite{gm} that if $F$ has characteristic zero then $L$ is completely factorisable if and only if the Frattini subalgebra of every subalgebra of $L$ is trivial. We shall see that this is false.
\par
If $A$ and $B$ are subalgebras of $L$ for which $L = A + B$ and $A \cap B = 0$ we will write $L = A \dot{+} B$; if, furthermore, $A, B$ are ideals of $L$ we write $L = A \oplus B$. The notation $A \leq B$ will indicate that $A$ is a subalgebra of $B$, and $A < B$ will mean that $A$ is a proper subalgebra of $B$. 
\bigskip
\section{Preliminary results}
\medskip
First we give some basic properties of c-supplemented subalgebras
\bigskip

\begin{lemma}\label{l:supp}
\begin{itemize}
\item[(i)] If $B$ is c-supplemented in $L$ and $B \leq K \leq L$ then $B$ is c-supplemented in $K$.
\item[(ii)] If $I$ is an ideal of $L$ and $I \leq B$ then $B$ is c-supplemented in $L$ if and only if $B/I$ is c-supplemented in $L/I$.
\item[(iii)] If ${\cal X}$ is the class of all c-supplemented Lie algebras then ${\cal X}$ is subalgebra and factor algebra closed.
\end{itemize}
\end{lemma}
\medskip
{\it Proof.}
\begin{itemize}
\item[(i)] Suppose that $B$ is c-supplemented in $L$ and $B \leq K \leq L$. Then there is a subalgebra $C$ of $L$ with $L = B + C$ and $B \cap C \leq B_L$. It follows that $K = (B + C) \cap K = B + C \cap K$ and $B \cap C \cap K \leq B_L \cap K \leq B_K$, and so $B$ is c-supplemented in $K$.
\item[(ii)] Suppose first that $B/I$ is c-supplemented in $L/I$. Then there is a subalgebra $C/I$ of $L/I$ such that $L/I = B/I + C/I$ and $(B/I) \cap (C/I) \leq (B/I)_{L/I} = B_L/I$. It follows that $L = B + C$ and $B \cap C \leq B_L$, whence $B$ is c-supplemented in $L$.
\par
Suppose conversely that $I$ is an ideal of $L$ with $I \leq B$ such that $B$ is c-supplemented in $L$. Then there is a subalgebra $C$ of $L$ such that $L = B + C$ and $B \cap C \leq B_L$. Now $L/I = B/I + (C + I)/I$ and
$(B/I) \cap (C + I)/I = (B \cap (C + I))/I = (I + B \cap C)/I \leq B_L/I = (B/I)_{L/I}$, and so $B/I$ is c-supplemented in $L/I$.
\item[(iii)] This follows immediately from (i) and (ii).
\end{itemize} 
\bigskip

The {\em Frattini ideal} of $L$, $\phi(L)$, is the largest ideal of $L$ contained in all maximal subalgebras of $L$. We say that $L$ is {\em $\phi$-free} if $\phi(L) = 0$. The next result shows that subalgebras of the Frattini ideal of a c-supplemented Lie algebra $L$ are necessarily ideals of $L$.
\bigskip

\begin{propo}\label{p:frat}
Let $B, D$ be subalgebras of $L$ with $B \leq \phi(D)$. If $B$ is c-supplemented in $L$ then $B$ is an ideal of $L$ and $B \leq \phi(L)$.
\end{propo}
\medskip
{\it Proof.} Suppose that $L = B + C$ and $B \cap C \leq B_L$. Then $D = D \cap L = D \cap (B + C) = B + D \cap C = D \cap C$ since $B \leq \phi(D)$. Hence $B \leq D \leq C$, giving $B = B \cap C \leq B_L$ and $B$ is an ideal of $L$. It then follows from \cite[Lemma 4.1]{frat} that $B \leq \phi(L)$.
\bigskip

The Lie algebra $L$ is called {\em elementary} if $\phi(B) = 0$ for every subalgebra $B$ of $L$; it is an {\em $E$-algebra} if $\phi(B) \leq \phi(L)$ for all subalgebras $B$ of $L$. Then we have the following useful corollary.
\bigskip

\begin{coro}\label{c:E}
If $L$ is c-supplemented then $L$ is an $E$-algebra.
\end{coro}
\medskip
{\it Proof.} Simply put $B = \phi(D)$ in Proposition \ref{p:frat}. 
\bigskip

It is clear that if $L$ is completely factorisable then it is c-supplemented. However, the converse is false. Every completely factorisable Lie algebra must be $\phi$-free, whereas the same is not true for c-supplemented algebras. For example, the three-dimensional Heisenberg algebra is c-supplemented, as will be clear from the next result which gives the true relationship between these two classes of algebras.
\bigskip

\begin{propo}\label{p:equ}
Let $L$ be a Lie algebra. Then the following are equivalent:
\begin{itemize}
\item[(i)] $L$ is c-supplemented.
\item[(ii)] $L/ \phi(L)$ is completely factorisable and every subalgebra of $\phi(L)$ is an ideal of $L$.
\end{itemize}
\end{propo}
\medskip
{\it Proof.} (i) $\Rightarrow$ (ii): Suppose first that $L$ is $\phi$-free and c-supplemented, and let $B$ be a subalgebra of $L$. Then there is a subalgebra $C$ of $L$ such that $L = B + C$. Choose $D$ to be a subalgebra of $L$ minimal with respect to $L = B + D$. Then $B \cap D \leq \phi(D)$, by \cite[Lemma 7.1]{frat}, whence $B \cap D = 0$ since $L$ is elementary, by Corollary \ref{c:E}. Hence $L$ is completely factorisable, and (ii) follows from Lemma \ref{l:supp}(iii) and Proposition \ref{p:frat}.
\par
(ii) $\Rightarrow$ (i): Suppose that (ii) holds and let $B$ be a subalgebra of $L$. Then there is a subalgebra $C/ \phi(L)$ of $L/ \phi(L)$ such that $L/ \phi(L) = ((B + \phi(L))/ \phi(L)) + (C/ \phi(L))$ and $ 0 = ((B + \phi(L))/ \phi(L)) \cap (C/ \phi(L)) = (B \cap C + \phi(L))/ \phi(L)$. Hence $L = B + C$ and $B \cap C \leq \phi(L)$, so $B \cap C$ is an ideal of $L$ and $B \cap C \leq B_L$; that is, $L$ is c-supplemented. 
\bigskip

Note that if $L$ is the three-dimensional Heisenberg algebra, then condition (ii) in the above result holds, since $\phi(L) = L^2$ is one dimensional and $L/ \phi(L)$ is abelian. Finally we shall need the following result concerning direct sums of of completely factorisable Lie algebras.
\bigskip

\begin{lemma}\label{l:dsum} 
If $A$ and $B$ are completely factorisable, then so is $L = A \oplus B$.
\end{lemma}
\medskip
{\it Proof.} Suppose that $A, B$ are completely factorisable and put $L = A \oplus B$. Let $U$ be a subalgebra of $L$. If $A \leq U$, then $U = A \oplus (B \cap U)$. Since $B$ is completely factorisable there is a subalgebra $C$ of $B$ such that $B = B \cap U + C$ and $U \cap C = B \cap U \cap C = 0$. Hence $L = U \dot{+} C$. 
\par
Now $A \leq A + U$ so, by the above, there is a subalgebra $C$ of $B$ with $L = A + U + C$ and $(A + U) \cap C = 0$. Moreover, since $A$ is completely factorisable, there is a subalgebra $D$ of $A$ such that $A = A \cap U + D$ and $U \cap D = A \cap U \cap D = 0$.  It follows that $L = U + (D \oplus C)$ and $U \cap (D + C) \leq U \cap [(A + U) \cap (D + C)] = U \cap [D + (A + U) \cap C] = U \cap D = 0$. It follows that $L$ is completely factorisable.
\bigskip

Note that the corresponding result for c-supplemented Lie algebras is false. For, let $L_1 = Fx + Fy + Fz$ with $[x,y] = -[y,x] = y + z, [x,z] = -[z,x] = z$ and all others products equal to zero. Then it is straightforward to check that $\phi(L_1) = Fz$ and that $L_1$ is c-supplemented. Now take $L$ to be a direct sum of two copies of $L_1$: say, $L = A \oplus B$ where $A = Fx + Fy + Fz$, $B = Fa + Fb + Fc$, $[x,y] = -[y,x] = y + z, [x,z] = -[z,x] = z$, $[a,b] = -[b,a] = b + c, [a,c] = -[c,a] = c$ and all others products equal to zero. Suppose that $F(z+c)$ is c-supplemented in $L$. Then there is a subalgebra $M$ of $L$ with $L = F(z+c) + M$
and $F(z+c) \cap M \leq (F(z+c))_L$. If $z+c \notin M$ then $M$ is a maximal subalgebra of $L$, contradicting the fact that $z+c \in (\phi(A) \oplus \phi(B)) = \phi(L)$, by \cite[Theorem 4.8]{frat}. It follows that $z+c \in M$, whence $F(z+c)$ is an ideal of $L$. But $[x, z+c] = z \notin F(z+c)$, a contradiction. Thus $L$ is not c-supplemented in $L$.
\bigskip

\section{The structure theorems}
\medskip
We can now give the main structure theorems for c-supplemented Lie algebras. First we determine the solvable ones.
\bigskip

\begin{theor}\label{t:solv}
Let $L$ be a solvable Lie algebra. Then the following are equivalent:
\begin{itemize}
\item[(i)] $L$ is c-supplemented.
\item[(ii)] $L$ is supersolvable and every subalgebra of $\phi(L)$ is an ideal of $L$.
\end{itemize}
\end{theor}
\medskip
{\it Proof.} (i) $\Rightarrow$ (ii): We have that every subalgebra of $\phi(L)$ is an ideal of $L$ by Proposition \ref{p:equ}, so we have only to show that $L$ is supersolvable. Let $L$ be a minimal counter-example. Then all proper subalgebras and factor algebras of $L$ are supersolvable, by Lemma \ref{l:supp}(iii). If we can show that all maximal subalgebras have codimension one in $L$, we shall have the desired contradiction, by \cite[Theorem 7]{barnes}; so let $M$ be any maximal subalgebra of $L$. Since the result is clear if $M_L \neq 0$, we may assume that $M_L = 0$.
\par
Pick a minimal ideal $A$ of $L$. Then $L = A \dot{+} M$ and $A$ is the unique minimal ideal of $L$, by \cite[Lemma 1.4]{tow}. Let $a \in A$. Then $Fa$ is c-supplemented in $L$, and so there is a subalgebra $B$ of $L$ such that $L = Fa + B$ and $Fa \cap B \leq (Fa)_L$. If $a \in B$ then $Fa$ is an ideal of $L$, whence $A = Fa$ and $M$ has codimension one in $L$.
\par
So suppose that $L = Fa \dot{+} B$. Since $A \not \leq B$ we have $B_L = 0$. But then $L = A \dot{+} B$ by \cite[Lemma 1.4]{tow} again. It follows that dim $A = 1$ and $M$ has codimension one in $L$.
\par
(ii) $\Rightarrow$ (i): By Proposition \ref{p:equ}, it suffices to show that if $L$ is supersolvable and $\phi$-free then it is completely factorisable. Let $L$ be a minimal counter-example. Then $L$ is elementary, by \cite[Theorem 1]{stit}, and so every proper subalgebra of $L$ is completely factorisable. Also $L = A \dot{+} B$ where $A = Fa_1 \oplus \ldots \oplus Fa_n$ is the abelian socle of $L$ and $B$ is abelian, by \cite[Theorem 7.3]{tow}. Let $U$ be a subalgebra of $L$. If $A \leq U$ it is clear that there is a subalgebra $C$ of $L$ such that $L = U + C$ and $U \cap C = 0$. So suppose that $a_i \notin U$ for some $1 \leq i \leq n$; we may as well assume that $i = 1$. Then $L/Fa_1 \cong (Fa_2 \oplus \ldots \oplus Fa_n) \dot{+} B$, which is a proper subalgebra of $L$ and so is completely factorisable. Hence there is a subalgebra $C$ of $L$ such that $L/Fa_1 = ((U + Fa_1)/Fa_1) + (C/Fa_1)$ and $Fa_1 = (U + Fa_1) \cap C = U \cap C + Fa_1$. It follows that $L = U + C$ and $U \cap C \leq Fa_1$. But $a_1 \notin U \cap C$ so $U \cap C = 0$ and $L$ is completely factorisable, a contradiction.
\bigskip 

We shall need the following classification of Lie algebras with core-free subalgebras of codimension one which is given by Amayo in \cite{am}.
\bigskip

\begin{theor}\label{t:am} (\cite[Theorem 3.1]{am})
Let $L$ have a core-free subalgebra of codimension one. Then either (i) dim $L \leq 2$, or else (ii) $L \cong L_m(\Gamma)$ for some $m$ and $\Gamma$ satisfying certain conditions (see \cite{am} for details).
\end{theor}
\bigskip

We shall also need the following properties of $L_m(\Gamma)$ which are given by Amayo in \cite{am}.
\bigskip

\begin{theor}\label{t:gamma} (\cite[Theorem 3.2]{am})
\begin{itemize}
\item[(i)] If $m > 1$ and $m$ is odd, then $L_m(\Gamma)$ is simple and has only one subalgebra of codimension one.
\item[(ii)] If $m > 1$ and $m$ is even, then $L_m(\Gamma)$ has a unique proper ideal of codimension one, which is simple, and precisely one other subalgebra of codimension one.
\item[(iii)] $L_1(\Gamma)$ has a basis $\{u_{-1}, u_0, u_1 \}$ with multiplication $[u_{-1}, u_0] = u_{-1} + \gamma_0 u_1$ $(\gamma_0 \in F, \gamma_0 = 0$ if $\Gamma = \{0\})$, $[u_{-1}, u_1] = u_0, [u_0, u_1] = u_1$.
\item[(iv)] If $F$ has characteristic different from two then $L_1(\Gamma) \cong L_1(0) \cong sl_2(F)$.
\item[(v)] If $F$ has characteristic two then $L_1(\Gamma) \cong L_1(0)$ if and only if $\gamma_0$ is a square in $F$. 
\end{itemize}
\end{theor}
\bigskip

The above properties enable us to determine which of the algebras $L_m(\Gamma)$ are c-supplemented.
\bigskip

\begin{propo}\label{p:gamma} 
If $L \cong L_m(\Gamma)$ then $L$ is c-supplemented if and only $L \cong L_1(0)$ and $F$ has characteristic different from two.
\end{propo}
\medskip
{\it Proof.} Suppose that $L \cong L_m(\Gamma)$ and $L$ is c-supplemented, and let $x \in L$. Then there is a subalgebra $M_1$ of $L$ such that $L = Fx + M_1$, and $Fx \cap M_1 \leq (Fx)_L = 0$, since $L_m(\Gamma)$ has no one-dimensional ideals. Choose $y \in M_1$. Then, similarly, there is a subalgebra $M_2$ of codimension one in $L$ such that $L = Fy + M_2$ and $M_1 \neq M_2$. Since $L = M_1 + M_2$ we have that $M_1 \cap M_2 \neq 0$. Let $z \in M_1 \cap M_2$. Then there is a subalgebra $M_3$ of codimension one in $L$ such that $L = Fz + M_3$, so $L$ has at least three subalgebras of codimension one in $L$. It follows from Theorem \ref{t:gamma} that $m = 1$.
\par
Suppose that $L \not \cong L_1(0)$. Then $F$ has characteristic two and $\gamma_0$ is not a square in $F$. Since $L$ is completely factorisable there is a two-dimensional subalgebra $M$ of $L$ such that $L = Fu_1 + M$. It follows that $M = F(u_{-1} + \alpha u_1) + F(u_0 + \beta u_1)$ for some $\alpha, \beta \in F$. But then $[u_{-1} + \alpha u_1, u_0 + \beta u_1] \in M$ shows that $\gamma_0 = \beta^2$, a contradiction. A further straightforward calculation shows that if $L \cong L_1(0)$ and $F$ has characteristic two, then $Fu_1$ is contained in every maximal subalgebra of $L$, and so has no c-supplement in $L$.
\par
Conversely, suppose that $L \cong L_1(0)$ and $F$ has characteristic different from two. Then $L \cong sl_2(F)$, by Theorem \ref{t:gamma} (iv) and it is easy to check that $L$ is c-supplemented.
\bigskip

We can now determine the simple and semisimple c-supplemented Lie algebras.
\bigskip

\begin{coro}\label{c:simple} 
If $L$ is simple then $L$ is c-supplemented if and only $L \cong L_1(0)$ and $F$ has characteristic different from two.
\end{coro}
\medskip
{\it Proof.} Let $L$ be simple and c-supplemented. Then $L$ has a core-free maximal subalgebra of codimension one in $L$ and so $L \cong L_m(\Gamma)$, by Theorem \ref{t:am}. The result now follows from Proposition \ref{p:gamma}.  
\bigskip

Notice, in particular, that $sl_2(F)$ is the only simple completely factorisable Lie algebra over any field. However, this is not the only simple elementary Lie algebra, even over a field of characteristic zero: over the real field every compact simple Lie algebra, and $so(n,1)$ for $n > 3$, for example, are elementary, as is shown in \cite[Theorem 5.1]{elem}. This justifies the assertion made at the end of the third paragraph of the introduction.
\bigskip

\begin{propo}\label{p:ss}
Let $L$ be a semisimple Lie algebra over a field $F$. Then the following are equivalent:
\begin{itemize}
\item[(i)] $L$ is c-supplemented.
\item[(ii)] $L = S_1 \oplus \ldots \oplus S_n$ where $S_i \cong sl_2(F)$ for $1 \leq i \leq n$ and $F$ has characteristic different from two.
\end{itemize}
\end{propo}
\medskip 
{\it Proof.} (i) $\Rightarrow$ (ii): Let $L$ be semisimple and c-supplemented and suppose the result holds for all such algebras of dimension less than dim $L$. Then $\phi(L) = 0$, since $\phi(L)$ is nilpotent, and so $L$ is completely factorisable. Let $A$ be a minimal ideal of $L$ and pick $a \in A$. Let $M$ be a subalgebra of $L$ such that $L = Fa \dot{+} M$ and put $B = A + M_L$. Then $M_L < B$ and $A \cap M_L = 0$, since $a \notin M_L$. If dim $L/M_L \leq 2$ then $A$ is abelian, contradicting the fact that $L$ is semisimple. It follows from Theorem \ref{t:am} and Proposition \ref{p:gamma} that $L/M_L \cong L_1(0)$, whence $B = L$ and $L = A \oplus M_L$. Since $A, M_L$ are semisimple and c-supplemented the result follows.
\par
(ii) $\Rightarrow$ (i):The converse follows from Corollary \ref{c:simple} and Lemma \ref{l:dsum}.  
\bigskip

Finally we have the main classification theorem.
\bigskip

\begin{theor}\label{t:supp}
Let $L$ be Lie algebra. Then the following are equivalent:
\begin{itemize}
\item[(i)] $L$ is c-supplemented.
\item[(ii)] $L/\phi(L) = R \oplus S$ where $R$ is supersolvable and $\phi$-free, $S$ is given by Proposition \ref{p:ss}, and every subalgebra of $\phi(L)$ is an ideal of $L$.
\end{itemize}
\end{theor}
\medskip
{\it Proof.} (i) $\Rightarrow$ (ii): Factor out $\phi(L)$ so that $L$ is $\phi$-free and c-supplemented and hence completely factorisable, by Proposition \ref{p:equ}. Then $L = R \dot{+} S$ where $R$ is the radical of $L$ and $S$ is semisimple. It suffices to show that $SR = 0$; the rest follows from Lemma \ref{l:supp}, Corollary \ref{c:E}, Proposition \ref{p:equ}, Theorem \ref{t:solv} and Proposition \ref{p:ss}. Suppose there is $0 \neq x \in L^{(3)} \cap R$. Then there is a subalgebra $M$ of $L$ such that $L = Fx \dot{+} M$ and $L/M_L$ is given by Theorem \ref{t:am}. If $L/M_L \cong L_m(\Gamma)$ then $L/ M_L$ is simple, by Proposition \ref{p:gamma}, and $M_L < R + M_L$, so $L = R + M_L$. But then $L/M_L$ is solvable, a contradiction. It follows that dim $L/M_L \leq 2$, whence $x \in L^{(3)} \cap R \leq L^{(3)} \leq M_L \leq M$, a contradiction. Hence $L^{(3)} \cap R = 0$. But $SR = S^2 R \leq S(SR) = S^2 (SR) \leq L^{(3)} \cap R = 0$, as required.
\par
(ii) $\Rightarrow$ (i): This follows from Proposition \ref{p:equ}, Lemma \ref{l:dsum}, Theorem \ref{t:solv} and Proposition \ref{p:ss}.

\end{document}